\documentclass{article}
\usepackage{amsmath}
\usepackage{amsthm}
\usepackage{amssymb}
\usepackage{mathrsfs}
\usepackage{amscd}
\usepackage{graphicx}
\usepackage{tikz}
\usepackage{indentfirst}

\newtheorem{thm}{Theorem}[section]
\newtheorem{cor}[thm]{Corollary}
\newtheorem{lem}[thm]{Lemma}

\newtheorem{claim}{Claim}

\begin{document}
\title{Ramsey numbers of color critical graphs versus large generalized fans\thanks{This work was supported by the National Natural Science Foundation of China (No. 12071453), the National Key R and D Program of China (2020YFA0713100),   and the Innovation Program for Quantum Science and Technology (2021ZD0302904).}}
\author{Taiping Jiang$^a$,\quad Xinmin Hou$^{a,b,c}$\\
\small $^a$School of Mathematical Sciences\\
\small University of Science and Technology of China, Hefei, Anhui 230026, China.\\
\small $^{b}$ CAS Key Laboratory of Wu Wen-Tsun Mathematics\\
\small University of Science and Technology of China, Hefei, Anhui 230026, China.\\
\small $^d$Hefei National Laboratory\\
\small University of Science and Technology of China, Hefei 230088, Anhui, China.
}

\date{}
\maketitle

\begin{abstract}
Given two graphs $G$ and $H$, the {Ramsey number} $R(G,H)$ is the smallest positive integer $N$ such that every 2-coloring of the edges of $K_{N}$ contains either a red $G$ or a blue $H$. Let $K_{N-1}\sqcup K_{1,k}$ be the graph obtained from $K_{N-1}$ by adding a new vertex $v$ connecting $k$ vertices of  $K_{N-1}$. Hook and Isaak (2011) defined the {\em star-critical Ramsey number} $r_{*}(G,H)$ as the smallest integer $k$ such that every 2-coloring of the edges of $K_{N-1}\sqcup K_{1,k}$ contains either a red $G$ or a blue $H$, where $N=R(G, H)$. 
For sufficiently large $n$, Li and Rousseau~(1996) proved that $R(K_{k+1},K_{1}+nK_{t})=knt +1$, Hao, Lin~(2018) showed that $r_{*}(K_{k+1},K_{1}+nK_{t})=(k-1)tn+t$;
 Li and Liu~(2016) proved that $R(C_{2k+1}, K_{1}+nK_{t})=2nt+1$, and Li, Li, and Wang~(2020) showed that $r_{*}(C_{2m+1},K_{1}+nK_{t})=nt+t$.  A graph $G$ with $\chi(G)=k+1$  is called edge-critical if  $G$ contains an edge $e$ such that $\chi(G-e)=k$.  In this paper, we extend the above results by showing that for an edge-critical graph $G$ with $\chi(G)=k+1$, when $k\geq 2$, $t\geq 2$ and $n$ is sufficiently large, $R(G, K_{1}+nK_{t})=knt+1$ and  $r_{*}(G,K_{1}+nK_{t})=(k-1)nt+t$. 

\end{abstract}

\section{Introduction}
Given two graphs $G$ and $H$, the {\em Ramsey number} $R(G,H)$ is the smallest positive integer $N$ such that every 2-coloring of the edges of $K_{N}$ contains either a red $G$ or a blue $H$.
Let $\chi(G)$ be the chromatic number of $G$, and $s(G)$ be the minimum size of a color class over all proper $\chi(G)$-coloring of $G$. Burr~\cite{ref2} gave  the following fundamental result about the lower bound of the Ramsey number $R(G, H)$.
\begin{thm}[Burr~\cite{ref2}]\label{THM: Burr}
Let $G$ and $H$ be graphs such that $H$ is connected and $|V(H)|\geq s(G)$. Then
\begin{equation}\label{EQ: lower}
R(G,H)\geq (\chi(G)-1)(|V(H)|-1)+s(G).
\end{equation}
\end{thm}

We call $H$ {\em $G$-good} if the lower bound in (\ref{EQ: lower}) is tight for $G$ and $H$. There are fruitful results on
Ramsey numbers $R(G, H)$, one can refer to an excellent  survey~\cite{ref12}. 
In this note, we are mainly concerned with  the pairs $(G, H)$ with $H$ being $G$-good (we also call  $(G, H)$ a {\em good pair}). Before we introduce  more results about good $(G, H)$-pairs, we give some notations. Write $K_k$ for a complete graph on $k$ vertices, and $T_n$ for a tree on $n$ vertices.  For disjoint graphs $G$ and $H$, let $G+H$ be the join graph of $G$ and $H$, i.e., the graph obtained by connecting $V(G)$ and $V(H)$ completely, and write  $nG$ for the disjoint union of $n$ copies of $G$. For example, a fan $F_{n}=K_{1}+nK_{2}$.  Write $C_k$ for a cycle on $k$ vertices.
There are some good pairs $(G, H)$ that we have known, here is a list of them.


\begin{itemize}
	\item[(a)]
	(Chv\'atal~\cite{ref3}) 
	Let $k,n\geq1$ be integers. Then $R(K_{k},T_{n})=(k-1)(n-1)+1$.
	
	\item[(b)] \textup{(Li and Rousseau~\cite{ref9})}
	If $k\geq 2$ and $t\geq 2$ are fixed, then $R(K_{k+1},K_{1}+nK_{t})=knt +1$ for sufficiently large $n$.
	
	\item[(c)] (Hamm, Hazelton, Thompson~\cite{HHS21}) Let $H$ be a graph on $h$ vertices, $k\ge 2$, $r\ge 1$. Then $R(rK_{k+1}, K_1+nH) = knh+r$  under
	certain restrictions on the relevant parameters.
	
	\item[(d)] (Nikiforov and Rousseau~\cite{NR04})
	If $k\geq 2$ and $t\geq 3$ are fixed, then $R(K_{k+1},K_{t}+nK_{1})=k(n+t-1) +1$ for sufficiently large $n$.
	
	\item[(e)] 
	\textup{(Li and Liu~\cite{ref10})}
	If integers $k,t\geq 1$ are fixed, then $R(C_{2k+1},K_{1}+nK_{t})=2nt+1$  for sufficiently large $n$.
	
\end{itemize}





 An interesting problem proposed by Hook and Isaak\cite{ref8} is the following: if we know the Ramsey number $R(G,H)=r$, what is the largest number of  edges we can delete from $K_{r}$ such that every 2-coloring of the edges of the resulting graph still contains either a red $G$ or a blue $H$. More precisely, let $K_{r-1}\sqcup K_{1,k}$ be the graph obtained from $K_{r-1}$ by adding a new vertex $v$ connecting $k$ vertices of  $K_{r-1}$.
Hook and Isaak~\cite{ref8} defined the {\em star-critical Ramsey number} $r_{*}(G,H)$ as the smallest integer $k$ such that every 2-coloring of the edges of $K_{r-1}\sqcup K_{1,k}$ contains either a red $G$ or a blue $H$. 
  There are a few known results of the star-critical Ramsey number, and we list some of them.
\begin{itemize}
\item[(f)] (Hao, Lin~\cite{ref7})
For any fixed integers $k\geq 2$ and $t\geq 2$, $r_{*}(K_{k+1},K_{1}+nK_{t})=(k-1)tn+t$ for all sufficiently large $n$.

\item[(g)] (Hamm, Hazelton, Thompson~\cite{HHS21}) Let $H$ be a graph on $h$ vertices, $k\ge 2$, $r\ge 1$. Then $r_*(rK_{k+1}, K_1+nH) = (k-1)nh+\delta(H)$  under
certain restrictions on the relevant parameters.

\item[(h)] (Hao, Lin~\cite{ref7})
For any fixed integers $k\geq 2$ and $t\geq 2$, $r_{*}(K_{k+1},K_{t}+nK_{1})=(k-1+o(1))n$ as $n\rightarrow\infty$. They also asked the question of the exact value of $r_{*}(K_{k+1},K_{t}+nK_{1})$.

\item[(i)](Li, Li, and Wang~\cite{ref11})
For a fixed odd cycle $C_{2m+1}$, if $n$ is large, then $r_{*}(C_{2m+1},K_{1}+nK_{t})=nt+t$.

\end{itemize}


A graph $H$ with $\chi(H)=k+1$ is called an edge-critical graph if $H$ contains an edge $e$ such that $\chi(H-e)=k$. Note that $C_{2m+1}$ and $K_{k+1}$ are edge-critical graphs. Recall that a fan $F_n=K_1+nK_2$. We call $K_1+nK_t$ a generalized fan. In this paper, we give the  Ramsey number and star-critical Ramsey number of an edge-critical graph versus a generalized fan. The following are the main results in this article.

\begin{thm}\label{THM: main1}
	Let $k\geq 2$, $t\geq 2$ be integers and $G$ be a fixed edge-critical graph with $\chi(G)=k+1$. If $n$ is sufficiently large, then $R(G,K_{1}+nK_{t})=knt+1$.
\end{thm}

\begin{thm}\label{THM: main2}
	Let $k\geq 2$, $t\geq 2$ be integers and $G$ be a fixed edge-critical graph with $\chi(G)=k+1$. If  $n$ is sufficiently large,
	then $r_{*}(G,K_{1}+nK_{t})=(k-1)nt+t$.
\end{thm}

Clearly, Theorems~\ref{THM: main1} and \ref{THM: main2}  extend the results (b), (e), (f), and (i). 

The rest of this article is arranged as follows. We provide some preliminaries in the next section. In Section 3, we prove Theorem~\ref{THM: main1} and give the proof of Theorem~\ref{THM: main2} in Section 4. We give some conclusions in the last section.


\section{Preliminaries}
In this section we will present some useful lemmas and results needed in our proofs.
The first lemma is an upper bound of $R(G,nH)$ due to Erd\H{o}s~\cite{ref1}.
\begin{lem}\textup{(\cite{ref1})}\label{LEM: upperbound}
Let $G$ and $H$ be graphs. Then $R(G,nH)\leq(n-1)|V(H)|+R(G,H)$.
\end{lem}

The following is the classical stability lemma due to Erd\H{o}s and Simonovits~\cite{ref4,ref5,ref6,ref13}. Let $G$ be a graph and $V_i, V_j\subseteq V(G)$. Write $E_G(V_i, V_j)$ for the set of edges with one end in $V_i$ and the other in $V_j$, and write $E_G(V_i)$ for $E_G(V_i, V_i)$. For a vertex $v\in V(G)$, write $d_{G}(v,V_i)=|E_G(\{v\}, V_i)|$. Let $[k]=\{1,2,\ldots,k\}$ for positive integer $k$.
\begin{lem}[Erd\H{o}s and Simonovits, 1966]\label{LEM: stability}
Let $G$ be a given graph with $\chi(G)=p+1$. For every $\xi>0$ there exist $\delta = \delta(\xi)$ and
$n_{0}=n_{0}(\delta)>0$ such that if $H$ is a graph of order $N>n_{0}$ and $e(H)>\dfrac{1}{2p}(p-1)N^{2}-\delta N^{2}$ that does not contain $G$, then there exists a partition of $V(H)$ into classes $V_{1}, V_{2},\ldots,V_{p}$ such that
\begin{itemize}
\item $\sum\limits_{1\le i<j\le p}\left(|V_i||V_j|-|E_H(V_i, V_j)|\right)<\xi N^{2}$;
\item $\sum\limits_{i=1}^p|E_H(V_i)|<\xi N^{2}$ for $i\in[p]$;
\item for every $v\in V_i$ and $i,j\in [p]$, $d_{G}(v, V_i)\le d_{G}(v, V_j)$;
\item $\frac Np - \xi N<|V_{i}|<\frac Np + \xi N$ for each $i\in[p]$.
\end{itemize}
\end{lem}

The following lemma obtained by Hao and Lin~\cite{ref7} gives a general lower bound for the star-critical Ramsey number $r_{*}(G,H)$.  
Denote a proper $k$-coloring of $G$ with color classes $U_{1}, U_{2},\ldots,U_{k}$ by $(U_{1}, U_{2},\ldots,U_{k})$. Let 
$$\tau(G)=\min_{v\in U_\chi\atop 1\le i<\chi}\{d_{G}(v, U_i) : \text{$(U_{1}, U_{2},\ldots,U_{\chi})$ is a proper $\chi$-coloring of $G$ with $|U_\chi|=s(G)$}\}.$$
\begin{lem}[\cite{ref7}]\label{LEM: 2.1r_*(G,H)}
		Let $G$ be a graph with $\chi(G)\geq 2$, and let $H$ be a connected graph of order $n\geq s(G)+1$ with minimum degree $\delta(H)$. We have
	$$r_{*}(G,H)\geq (\chi(G)-2)(n-1)+\min\{n,\delta(H)+\tau(G)-1\}.$$
Furthermore, if $H$ contains no  cut vertex or $\delta(H)=1$, then
	$$r_{*}(G,H)\geq (\chi(G)-2)(n-1)+\min\{n,\delta(H)+\tau(G)-1\}+s(G)+1.$$
	
\end{lem}
We also call $(G, H)$ a good pair if the lower bound in the above lemma is tight for $r_*(G, H)$.
Finally, we give a structural property of the edge-critical graph. 
\begin{lem}\label{LEM: 2.4criticalstructrue}
Let $G$ be an edge-critical graph with $\chi(G)=k$. Then there exists a proper $k$-coloring $(V_1, V_2, \ldots, V_{k})$ of $G$ such that $|V_k|=1$ and $|E(V_k, V_i)|=1$ for some $i\in[k-1]$. 
\end{lem}
\begin{proof}

Choose an edge $e=uv\in E(G)$ with $\chi(G-e)=k-1$. Then $G-e$ has a proper $(k-1)$-coloring $\pi'$ with $\pi'(u)=\pi'(v)$, otherwise, we will obtain a proper $(k-1)$-coloring of $G$,  a contradiction. 
Now we recolor the vertices by setting $\pi:  V(G)\rightarrow [k]$ as follows:
\begin{equation}
\pi(x)=\left\{
\begin{array}{rr}
\pi'(x), & {x\neq v},\nonumber \\
k,        & {x=v}.
\end{array}
\right.
\end{equation}
Then $\pi$ is a proper $k$-coloring of $G$, and the vertex $v$ is what we want to find.
\end{proof}

As a direct corollary of Lemma~\ref{LEM: 2.4criticalstructrue}, we have 
\begin{cor}\label{COR: s(G)=1}
Let $G$ be an edge-critical graph. Then $s(G)=1$ and $\tau(G)=1$.
\end{cor}



\section{Proof of Theorem~\ref{THM: main1}: $R(G,K_{1}+nK_{t})=knt+1$}
Write $K_{n_1,n_2,\ldots, n_k}$ for a $k$-partite graph with color classes of sizes $n_1, n_2,\ldots, n_k$ and let $K^+_{n_1,n_2,\ldots, n_k}$ be the graph obtained by adding an extra edge in one class of $K_{n_1,n_2,\ldots, n_k}$.

\noindent
\textbf{Proof of Theorem 3.1.} 
Let $G$ be an edge-critical graph with $\chi(G)=k+1$. By Corollary~\ref{COR: s(G)=1}, $s(G)=1$.
By Theorem~\ref{THM: Burr}, $R(G,K_{1}+nK_t)\geq knt+1$. Let $N=knt+1$. Then it suffices to show $R(G,K_{1}+nK_t)\leq knt+1$, i.e., we will prove that any 2-coloring of the edges of $K_{N}$ contains either a red $G$ or a blue $K_{1}+nK_{t}$ when $n$ is sufficiently large.

Now suppose to the contrary that there exists a 2-coloring of the edges of $K_N$ that contains neither a red $G$ nor a blue $K_{1}+nK_{t}$. Let $R$ and $B$ denote the subgraphs consisting of red and blue edges in $K_N$, respectively.   
Then the subgraph of $K_N$ induced by $N_B(v)$ contains neither a red $G$ nor a blue $nK_{t}$. By Lemma~\ref{LEM: upperbound}, we have
\begin{center}
$d_{B}(v)<R(G,nK_{t})\leq (n-1)t+R(G,K_t)$.
\end{center}
Thus $d_{R}(v)=N-1-d_{B}(v)\geq (k-1)tn-O(1)$, which implies that
\begin{center}
$e(R)=\dfrac{1}{2}\sum\limits_{v}d_{R}(v)\geq \dfrac{k-1}{2k}N^{2}-O(N)$.
\end{center}

Set $\xi>0$, a very small real number. Applying Lemma~\ref{LEM: stability} to graph $G$ and the $G$-free graph $R$, we have 
a partition of $V(R)= V(K_{N})$ into $k$ classes $V_{1},V_{2},\cdots,V_{k}$ such that
\begin{itemize}
\item[(i)] $\sum\limits_{1\le i<j\le k}\left(|V_i||V_j|-|E_R(V_i, V_j)|\right)<\xi N^{2}$;
\item[(ii)] $|E_R(V_i)|<\xi N^{2}$ for every $i\in[k]$;
\item[(iii)] for every $v\in V_i$ and $i,j\in [k]$, $d_{R}(v, V_i)\le d_{R}(v, V_j)$;
\item[(iv)] $\frac Nk-\xi N<|V_{i}|<\frac Nk + \xi N$ for every $i\in [k]$.
\end{itemize}
With loss of generality, we may assume that $|V_{1}|\geq |V_{i}|$ for $i=2,\cdots,k$. Thus $|V_{1}|\geq tn+1>\frac Nk$. Define a subset $V_{i}'\subset V_{i}$ as follows:
\begin{center}
$V_{i}'=\{x\in V_{i}\ |\ d_{R}(x,V_{j})\geq (1-2\sqrt{\xi})|V_{j}|,\text{ for all } j\in[k] \text{ and } j\neq i\}$.
\end{center}
\begin{claim}\label{Claim 1.}
	 $|V_{i}'|\geq(1-2k\sqrt{\xi})|V_{i}|$ for $i\in [k]$.
\end{claim}
\textbf{Proof.} Suppose to the contrary that there is an $i\in[k]$ with 
\begin{center}
$|V_{i}\setminus V_{i}'|>2k\sqrt{\xi}|V_{i}|>2k\sqrt{\xi}\left(\frac 1k-\xi\right)N$.
\end{center}
By the definition of $V_{i}'$, for any vertex $v\in V_{i}\setminus V_{i}'$, we have
\begin{center}
$d_{B}(v,V_{j})=|V_j|-d_R(x, V_j)\geq 2\sqrt{\xi}|V_{j}|\geq 2\sqrt{\xi}\left(\frac 1k-\xi\right)N$,  for all  $j\in[k]$  and  $j\neq i$.
\end{center}
Thus 
\begin{center}
$\sum\limits_{j=1\atop j\not=i}^k|E_B(V_i, V_j)|>(k-1)\cdot2k\sqrt{\xi}\left(\frac 1k-\xi\right)N\cdot2\sqrt{\xi}\left(\frac 1k-\xi\right)N=4k(k-1)\xi\left(\frac 1k-\xi\right)^{2}N^{2}>\xi N^{2}$.
\end{center}
However, $$\sum\limits_{j=1\atop j\not=i}^k|E_B(V_i, V_j)|\le \sum\limits_{1\le i<j\le k}|E_B(V_i, V_j)|=\sum\limits_{1\le i<j\le k}\left(|V_i||V_j|-|E_R(V_i, V_j)|\right)<\xi N^{2},$$
 a contradiction.   $\hfill\square$
  
Let $m$ be the maximum size of a color class over all proper $(k+1)$-coloring of $G$. The construction of a copy of red $G$ in $R$ will be used throughout the rest of the proof.
\begin{claim}\label{construction}
If there is a red edge $e=u_{i1}u_{i2}$	in some $V_i$ such that $d_R(u_{i1}, V_j')+d_R(u_{i2}, V_j')\ge (1+\alpha)|V_j|$ for some  constant $\alpha>0$ and $j\not=i$, then $R$ contains a copy of $G$.
\end{claim}	
\begin{proof}
Without loss of generality, we assume that $e=u_{11}u_{12}$ is a red edge in $V_1$.  Note that, for any vertex $v\in V_{i}'$, we have 
	\begin{center}
		$d_{R}(v,V_{j}')\geq (1-2\sqrt{\xi})|V_{j}|-2k\sqrt{\xi}|V_{j}|=(1-2(k+1)\sqrt{\xi})|V_{j}|$.
	\end{center}
	We will construct a copy of $G$ in $R$ with $e=u_{11}u_{12}\in E(G)$. In the following proof of this claim, the host graph we considered is $R$. First, we choose $m-1$ vertices $u_{13}, \ldots, u_{1(m+1)}$ from $V_1'$ and let $U_{1}=\{u_{11},u_{12},u_{13}\ldots,u_{1(m+1)}\}$. Then, when $n$ is large enough, the common neighbors of $U_1$ in $V_j'$ with $j\ge 2$ are 
	\begin{center}
		$\bigg\lvert\bigcap\limits_{i=1}^{m+1}N_{R}(u_{1i},V_{j}')\bigg\rvert\geq  (1+\alpha)|V_j|+ \sum\limits_{i=3}^{m+1}d_R(u_{1i}, V_j')-m|V_j'| \geq \left(\alpha-2(k+1)(m-1)\sqrt{\xi}\right)|V_{j}|\geq m$.
	\end{center}
	For ease of notations, we denote $\bigcap\limits_{i=1}^{m+1}N_{R}(u_{1i},V_{j}')=W_j^{(1)}$ for $j=2,3,\ldots,k$. Second, we take $m$ vertices 
	$U_{2}=\{u_{21},u_{22},\cdots,u_{2m}\}$ from $W_2^{(1)}$. Then, when $n$ is large enough, the common neighbors of $U_{1}\cup U_{2}$ in $V_j'$ for $j=3,\ldots, k$ is 
	\begin{center}
		$\bigg\lvert\left(\bigcap\limits_{i=1}^{m+1}N_{R}(u_{1i},V_{j}')\right)\bigcap\left( \bigcap\limits_{i=1}^{m}N_{R}(u_{2i},V_{j}')\right)\bigg\rvert\ \geq \left(\alpha-2(k+1)(2m-1)\sqrt{\xi}\right)|V_{j}|\geq m$.
	\end{center}
	Similiarly, denote $\left(\bigcap\limits_{i=1}^{m+1}N_{R}(u_{1i},V_{j}')\right)\bigcap \left(\bigcap\limits_{i=1}^{m}N_{R}(u_{2i},V_{j}')\right)=W_j^{(2)}$,
	for $j=3,\ldots,k$. Note that $W_3^{(2)}\subseteq W_3^{(1)}$. We continue to choose $m$ vertices 
	$U_{3}=\{u_{31},u_{32},\cdots,u_{3m}\}$ from $W_3^{(2)}$ and compute the common neighbors of $U_{1}\cup U_{2}\cup U_3$ in $V_{4}',\ldots,V_{k}'$. Repeat the above process, until the last step we get $k-1$  disjoint vertex sets $U_{1},U_{2},\cdots,U_{k-1}$, consider their common neighbors in $V_{k}'$.
	Let $$\left(\bigcap\limits_{i=1}^{m+1}N_{R}(u_{1i},V_{k}')\right)\bigcap \left(\bigcap\limits_{i=1}^{m}N_{R}(u_{2i},V_{k}')\right)\bigcap \cdots \bigcap\left(\bigcap\limits_{i=1}^{m}N_{R}(u_{(k-1)i},V_{k}')\right)=W_k^{(k-1)}.$$
	When $n$ is large enough,  we still have 
	\begin{center}
		$|W_k^{(k-1)}|\geq\left(\alpha-2(k+1)((k-1)m-1)\sqrt{\xi}\right)|V_{k}|\geq m$.
	\end{center}
	Thus we  can still choose $m$ vertices
	$U_{k}=\{u_{k1},u_{k2},\cdots,u_{km}\}$ from $W_k^{(k-1)}$. Note that $W_i^{(j)}\subseteq W_i^{(j')}$ for $2\le i\le k$ and $1\le j'\le j\le k-1$. Thus $R[U_{1}\cup U_{2}\cup\ldots\cup U_{k}]$ contains a copy of $K^+_{(m+1),m,\ldots,m}$ with an extra edge $e=u_{11}u_{12}$ in $R[U_1]$. By the definition of $m$ and Lemma~\ref{LEM: 2.4criticalstructrue}, we know that there is a copy of $G$ in $K^+_{(m+1),m,\ldots,m}\subseteq R$.
\end{proof}
 
\begin{claim}\label{Claim 2.} 
	The induced subgraph $B[V_{i}']$ is a complete graph for $i\in [k]$.
\end{claim}
\noindent\textbf{Proof.} By symmetry, it is sufficient to show that $B[V_{1}']$ is a complete graph. Suppose not, then there exists an edge $e=\{u_{11},u_{12}\}$ in $R[V_{1}']$. Note that, for any vertex $v\in V_{i}'$, we have 
\begin{center}
$d_{R}(v,V_{j}')\geq (1-2\sqrt{\xi})|V_{j}|-2k\sqrt{\xi}|V_{j}|=(1-2(k+1)\sqrt{\xi})|V_{j}|$.
\end{center}
Thus        $d_R(u_{11}, V_j')+d_R(u_{12}, V_j')\ge (2-4(k+1)\sqrt{\xi})|V_j|\ge (1+\alpha)|V_j|$ for some constant $\alpha>0$.
By Claim~\ref{construction}, we have a red  copy of $G$ in $R$, a contradiction.

\begin{claim}\label{Claim 3.} 
	$E(V_{i}',V_{i}\backslash V_{i}')\subseteq E(B)$ for every $i\in[k]$.
\end{claim}
\noindent\textbf{Proof.} By symmetry, it suffices to prove that $E(V_{1}', V_{1}\backslash V_{1}')\subseteq E(B)$. Suppose to the contrary that there exists a red edge $e=u_{11}u_{12}$ with $u_{11}\in V_{1}\backslash V_{1}'$ and $u_{12}\in V_{1}'$. 
Again note that,   for $j=2,3,\ldots, k$,
\begin{center}
	$d_{R}(u_{12},V_{j}')\geq (1-2\sqrt{\xi})|V_{j}|-2k\sqrt{\xi}|V_{j}|=(1-2(k+1)\sqrt{\xi})|V_{j}|$.
\end{center}
By Lemma~\ref{LEM: stability}, $d_{R}(u_{11},V_{j})\ge d_R(u_{11}, V_1)$ for $u_{11}\in V_1\setminus V_1'$. 
Thus, when $\xi$ is very small and $n$ is large enough, we have 
\begin{eqnarray*}
d_{R}(u_{11},V_{j}')&\geq& \dfrac{1}{2}\left(d_{R}(u_{11}, V_1)+d_{R}(u_{11},V_{j})\right)-|V_j\setminus V_j'|\\
              &=& \dfrac{1}{2}\left(d_R(u_{11})-\sum_{i=2\atop i\not=j}^k|V_i|\right)-|V_j\setminus V_j'|\\
              &\geq& \frac 12\left((k-1)nt-O(1)-(k-2)\left(\frac Nk+\xi N\right)\right)-|V_j\setminus V_j'|\\
              &\geq& \frac 12\left((k-1)\frac{N-1}k-O(1)-(k-2)\left(\frac Nk+\xi N\right)\right)-2k\sqrt\xi |V_j|\\
              &>&\left(\dfrac{1}{3}-2k\sqrt\xi\right)|V_{j}|.
\end{eqnarray*}
Therefore,  $d_R(u_{11}, V_j')+d_R(u_{12}, V_j')\ge (\frac 43-2(2k+1)\sqrt{\xi})|V_j|\ge (1+\alpha)|V_j|$ for some constant $\alpha>0$. Again by Lemma~\ref{construction}, we have a red copy of $G$ in $R$, a contradiction.

\vspace{5pt}

By Claims~\ref{Claim 2.} and~\ref{Claim 3.}, any vertex $v\in V_{1}\backslash V_{1}'$ and $t-1$ vertices in $V_{1}'$ induces  a complete graph $K_{t}$ in $B$. Since $|V_{1}$\textbackslash$ V_{1}'|<|V_{1}'|/(t-1)$, we can repeat the above process until the vertices of $V_{1}$\textbackslash$ V_{1}'$
have been used up, now we have $|V_1\setminus V_1'|$ copies of blue $K_t$ in $B[V_1]$.  Note that $|V_1|\ge nt+1$. By Claim~\ref{Claim 2.}, we can continue to find copies of blue $K_t$ in the remaining vertices of $B[V_{1}]$ until we find $n$ disjoint copies of blue $K_{t}$ in $B[V_1]$. Note that $|V_{1}|\geq nt+1$, there is at least one vertex left in $V_1'$ that is completely connected with the $n$ vertex disjoint copies of blue $K_{t}$, thus we construct a blue $K_{1}+nK_t$, which leads to a contradiction.
$\hfill\square$

\section{Proof of Theorem~\ref{THM: main2}: $r_{*}(G,K_{1}+nK_{t})=(k-1)nt+t$}

\noindent
\textbf{Proof of Theorem~\ref{THM: main2}:} By Corollary~\ref{COR: s(G)=1}, we have $\tau(G)=1$.  Note that $\delta(K_{1}+nK_{t})=t$. Thus when $n$ is large, we have 
$r_{*}(G,K_{1}+nK_{t})\geq(k-1)nt+t$ by Lemma~\ref{LEM: 2.1r_*(G,H)}. In the following, we will show that $r_{*}(G,K_{1}+nK_{t})\leq(k-1)nt+t$.

Let $N=knt$. Then $R(G,K_{1}+nK_{t})=N+1$ by Theorem~\ref{THM: main1}. Thus there is a $\{R, B\}$-coloring  of the edges of $K_N$ such that it contains neither a red $G$ nor a blue $K_{1}+nK_{t}$.
By Lemma~\ref{LEM: upperbound}, we have
\begin{center}
$d_{B}(v)<R(G,nK_{t})\leq t(n-1)+R(G,K_t)$.
\end{center}
Therefore, $d_{R}(v)=N-1-d_{B}(v)\geq (k-1)tn-O(1)$ which implies that
\begin{center}
$e(R)=\dfrac{1}{2}\sum\limits_{v\in V(R)}d_{R}(v)\geq \dfrac{k-1}{2k}N^{2}-O(N)$.
\end{center}
Apply the stability Lemma~\ref{LEM: stability} to $R$  with the forbidden graph $G$, for a sufficiently small real number $\xi>0$,
there is a partition of $V(R)= V(K_{N})$ into $k$ classes $V_{1},V_{2},\ldots,V_{k}$ such that
\begin{itemize}
	\item[(i)] $\sum\limits_{1\le i<j\le k}\left(|V_i||V_j|-|E_R(V_i, V_j)|\right)<\xi N^{2}$;
	\item[(ii)] $|E_R(V_i)|<\xi N^{2}$ for every $i\in[k]$;
	\item[(iii)] for every $v\in V_i$ and $i\in [k]$, $d_{R}(v, V_i)\le d_{R}(v, V_j)$ for every $j\in[k]$;
	\item[(iv)] $\frac Nk-\xi N<|V_{i}|<\frac Nk + \xi N$ for every $i\in [k]$.
\end{itemize}


Similar to the proof of Theorem~\ref{THM: main1}, define $V_{i}'\subset V_{i}$ as follows:
\begin{center}
$V_{i}'=\{x\in V_{i}\, |\, d_{R}(x,V_{j})\geq (1-2\sqrt{\xi})|V_{j}|,\text{ for all $j\in [k]$ and } j\neq i\}$.
\end{center}
With the same discussion as in the proof of Theorem~\ref{THM: main1}, we have the following three claims.

\begin{claim}\label{P2:Claim 1.} 
	$|V_{i}'|\geq(1-2k\sqrt{\xi})|V_{i}|$ for every $i\in[k]$.
\end{claim}

\begin{claim}\label{P2:Claim 2.}
 $B[V_{i}']$ is a complete graph for every $i\in[k]$.
\end{claim}

\begin{claim}\label{P2:Claim 3.} 
	$E(V_{i}', V_{i})\subseteq E(B)$ for every $i\in[k]$.
\end{claim}
By Claims~\ref{P2:Claim 2.} and~\ref{P2:Claim 3.}, we have $|V_{i}|\leq nt$, otherwise, we can construct a copy of $K_{1}+nK_{t}$ in $B$ with the same discussion as in the last paragraph of the proof of Theorem~\ref{THM: main1}. Note that $V_{1},V_{2},\ldots,V_{k}$ is a partition of $V({K_{N}})$ and $N=knt$. We have the following claim. 

\begin{claim}\label{CL: |V_i|=nt}
For every $i\in [k]$, $|V_{i}|=nt$ and there are $n$ disjoint copies of blue $K_t$ in $B[V_i]$.
\end{claim}
We further claim that there are no red edges in $V_i\setminus V_i'$ for $i\in[k]$. Recall that $m$ denotes the maximum size of a color class over all proper $(k+1)$-coloring of $G$.
\begin{claim}\label{P2:Claim 4.}
	 $B[V_{i}\backslash V_{i}']$ is a complete graph for every $i\in[k]$.
\end{claim}
\noindent\textbf{Proof.} By symmetry, we will show that $B[V_{1}\backslash V_{1}']$ is a complete graph. Suppose to 
the contrary that there exists a red edge $v_{11}v_{12}$ with $v_{11},v_{12}\in V_{1}$\textbackslash$ V_{1}'$. Then, $d_{B}(v_{11},V_{j}')\leq t-1$  for $j=2,\ldots, k$, since, otherwise, $\{v_{11}\}\cup V_j$ contains a copy of blue $K_{1}+nK_{t}$. 
Similarly, we have $d_{B}(v_{12},V_{j}')\leq t-1$. Thus, we have
\begin{center}
	$d_{R}(v_{1i},V_{j}')\ge |V_j'|-t+1\ge (1-2k\sqrt{\xi})|V_{j}|-t+1$ for $i=1,2$.
\end{center}
Therefore, when $n$ is sufficiently large and $\xi$ is very small, 
 $$d_{R}(v_{11},V_{j}')+d_R(v_{12}, V_j')\ge (2-4k\sqrt{\xi})|V_{j}|-2t+2\ge (1+\alpha)|V_j|$$ for some constant $\alpha>0$.
 By Claim~\ref{construction}, $R$ contains a copy of $G$, a contradiction.

By Claims~\ref{P2:Claim 2.}, \ref{P2:Claim 3.}, and \ref{P2:Claim 4.}, the graph induced by $V_{i}$ is a blue $K_{tn}$ for $i\in[k]$.
\vspace{5pt}

Now, assume $v_{0}$ is a new vertex outside $V({K_{N}})$. We claim that if $d(v_{0}, V(K_N))\ge (k-1)nt+t$, then there exists either a red $G$ or a blue $K_{1}+nK_{t}$.

Without loss of generality, suppose $d(v_{0},V_{k})\leq d(v_{0},V_{k-1})\leq \cdots\leq d(v_{0},V_{1})$. Then we have
$t\leq d(v_{0},V_{k})\leq (k-1)nt/k+t/k$. By Claim~\ref{CL: |V_i|=nt}, $d(v_0, V_i)\le |V_i|=nt$. We can easily prove that $d(v_{0},V_{k-1})\geq nt/2$, since, otherwise,
\begin{center}
$\sum\limits_{i=1}^{k-2}d(v_{0},V_{i})+d(v_{0},V_{k-1})+d(v_{0},V_{k})<(k-2)nt+nt/2+nt/2=(k-1)nt<d(v_{0}, V(K_N))$, 
\end{center}
a contradiction. Thus we have $d(v_{0},V_{i})\ge d(v_0, V_{k-1})\geq nt/2$ for $i\in[k-1]$. Note that $V_i$ induces a blue $K_{nt}$. Then for each $i\in[k]$, $d_{B}(v_{0},V_{i})\leq t-1$, otherwise, $\{v_0\}\cup V_i$ contains a blue $K_{1}+nK_{t}$. 
For the same reason, we know that, for $u\in V_i$, $d_B(u, V_j)\le t-1$ for every $j\in[k]$ and $j\not=i$.
Therefore, 
\begin{center}
$d_{R}(v_{0},V_{k})\geq d(v_0, V_k)-(t-1)\ge t-(t-1)=1$,
\end{center}
and $$ \text{ for  $u\in V_i$, }  d_R(u, V_j')\ge |V_j'|-(t-1)\ge (1-2k\sqrt{\xi})|V_j|-t+1, \text{ $j\in[k]$ and $j\not=i$.}$$
Thus there exists a vertex $v_{11}\in V_{k}$ such that $v_{0}v_{11}$ is a red edge. Note that 
$$d_{R}(v_{0},V_{j}')\geq \frac{nt}2-|V_j\setminus V_j'| -(t-1)\ge \left(\frac 12-2k\sqrt{\xi}\right)|V_j|-t+1, \text{ for } j\in [k-1].$$
Therefore, when $n$ is sufficiently large and $\xi$ is very small, we have 
$$d_{R}(v_{0},V_{j}')+d_R(v_{11}, V_j')\geq \left(\frac 32-4k\sqrt{\xi}\right)|V_j|-2t+2\ge(1+\alpha)|V_j|,$$
\text{ for  some constant $\alpha>0$ and } $j\in [k-1].$
We can view $v_0$ as a vertex in $V_k$ and $v_0v_{11}$ as a red edge in $V_k$. Then, by Claim~\ref{construction}, we have a copy of red $G$ in $R$, a contradiction.
This completes the proof of Theorem 3.2. $\hfill\square$

\section{Discussion and remarks}

In this paper, we show that for an edge-critical graph $G$ with $\chi(G)=k+1$, when $k\geq 2$, $t\geq 2$ and $n$ is sufficiently large, $R(G, K_{1}+nK_{t})=knt+1$ and  $r_{*}(G,K_{1}+nK_{t})=(k-1)nt+t$. 
As mentioned in the introduction, Nikiforov and Rousseau~(see (d)) proved that $R(K_{k+1}, K_{t}+nK_{1})=k(n+t-1)+1$ for sufficiently large $n$, and Hao and Lin~(see (h)) showed that  $r_{*}(K_{k+1},K_{t}+nK_{1})=(k-1+o(1))n$ as $n\rightarrow\infty$ and asked the question of what the exact value of $r_{*}(K_{k+1},K_{t}+nK_{1})$ is.
One may naturally ask the question of what the exact values of $(G, K_t+nK_1)$ and $r_{*}(G,K_{t}+nK_{1})$ are? More generally,  determine the exact values of $(G, K_t+nK_r)$ and $r_{*}(G,K_{t}+nK_{r})$.

\end{document}